\documentclass[10pt,notitlepage,twoside,a4paper]{amsart}
 \usepackage{amsfonts}

\usepackage{amsmath,amssymb,enumerate}

\usepackage{epsfig,fancyhdr,color}

\usepackage{amssymb}
\usepackage{amsmath,amsthm}
\usepackage{latexsym}
\usepackage{amscd,mathrsfs,marginnote}
\usepackage{psfrag}
\usepackage{graphicx}
\usepackage[latin1]{inputenc}
\usepackage[all]{xy}
\usepackage[mathcal]{eucal}
\usepackage[margin=2in]{geometry}

\definecolor{NoteColor}{rgb}{1,0,0}

\makeatletter
\newcommand\Max{\@tempcnta=\mathcode`\m\relax
\mathcode`\m=\mathcode`\M\max\mathcode`\m=\@tempcnta\relax}

\newcommand\Min{\@tempcnta=\mathcode`\m\relax
\mathcode`\m=\mathcode`\M\min\mathcode`\m=\@tempcnta\relax}
\makeatother

\newtheorem{theorem}{\rm\bf Theorem}[section]
\newtheorem{proposition}[theorem]{\rm\bf Proposition}

\newtheorem*{theorem 1}{\rm\bf Proposition 1}
\newtheorem*{theorem 2}{\rm\bf Proposition 2}

\theoremstyle{definition}

\theoremstyle{remark}

\newtheorem{remarks}[theorem]{\rm\bf Remarks}

\def\interieur#1{\mathord{\mathop{\kern 0pt #1}\limits^\circ}}

\makeindex

\title[The type problem]{Teichm\"uller's work on the type problem}

\author[V. Alberge, M. Brakalova and A. Papadopoulos]{Vincent Alberge, Melkana Brakalova-
\\
 and Athanase Papadopoulos}

\thanks{}
 

\date{\today}

\begin{document}

  \maketitle
  
  \begin{abstract}
  The \emph{type problem} is the problem of deciding, for a simply connected Riemann surface, whether it is conformally equivalent to the complex plane or to the unit dic in the complex plane.
  
 We report on Teichm\"uller's results on the type problem from his two papers \emph{Eine Anwendung quasikonformer Abbildungen auf das Typenproblem} (An application of quasiconformal mappings to the type problem) (1937) and \emph{Untersuchungen \"uber konforme und quasikonforme Abbildungen} (Investigations on conformal and quasiconformal mappings) (1938). They concern simply connected Riemann surfaces defined as branched covers of the sphere.
  At the same time, we review the theory of line complexes, a combinatorial device used by Teichm\"uller and others to encode branched coverings of the sphere. 
  
  In the first paper, Teichm\"uller proves that any two simply connected Riemann surfaces which are branched coverings of the Riemann sphere with finitely many branch values and which have the same line complex are quasiconformally equivalent. For this purpose, he introduces a technique for piecing together quasiconformal mappings. He also obtains a result on the extension of smooth orientation-preserving diffeomorphisms of the circle to quasiconformal mappings of the disc which are conformal at the boundary. 
  
  In the second paper, using line complexes, Teichm\"uller gives a type criterion for a simply-connected surface which is a branched covering of the sphere, in terms of an adequately defined measure of ramification, defined by a limiting process. The result says that if the surface is ``sufficiently ramified" (in a sense to be made precise), then it is hyperbolic. In the same paper, Teichm\"uller answers by the negative a conjecture made by Nevanlinna\index{Nevanlinna, Rolf (1895--1980)} which states a criterion for parabolicity in terms of the value of a (different) measure of ramification, defined by a limiting process. Teichm\"uller's results in his first paper are used in the proof of the results of the second one.

 The final version of this paper will appear in Vol. VII of the \emph{Handbook of Teichm\"uller theory} (European Mathematical Society Publishing House, 2020).

\end{abstract}

\bigskip
 
\bigskip

 \section{Introduction: The type problem}

   By the\index{type problem} uniformization theorem, every simply connected Riemann surface $S$ is either conformally equivalent to the unit disc, or to the complex plane, or to the Riemann sphere. In the first case, the surface is said to be of hyperbolic type, in the second, of parabolic type, and in the third, of elliptic type. The last case is distinguished from the two others by topology: if $S$ is contractible, it is not homeomorphic to the Riemann sphere. How can we distinguish between the other two cases? In other words, suppose that a contractible Riemann surface $S$ is  constructed by some process, for example, by assembling pieces which are themselves simply connected surfaces, or by analytic continuation, or---as in the papers considered---as an infinite branched covering of the sphere with information on the ramification index at the branching points. How can we decide if the surface is of hyperbolic or parabolic type?  This is the so-called type problem.

   Ahlfors\index{Ahlfors, Lars (1907--1996)} writes in his comments to his \emph{Collected works} \cite[vol. 1, p. 84]{Ahlfors-Collected}, that the type problem\index{type problem} was formulated for the first time by Andreas Speiser\footnote{Andreas Speiser\index{Andreas, Speiser (1885--1970)} (1885--1970) was a Swiss mathematician. He studied in G\"ottingen, first with Minkowski and then with Hilbert, who became his doctoral advisor after the death of Minkowski. He defended his PhD thesis, in Göttingen, in 1909 and his habilitation in 1911, in Strasbourg. Speiser\index{Speiser, Andreas (1885--1970)} worked in group theory, number theory and Riemann surfaces, but he is mostly known as the main editor of Euler's \emph{Opera Omnia}. (He edited 11 volumes and collaborated to 26 others.)} in his paper \emph{Probleme aus dem Gebiet der ganzen transzendenten Funktionen} (Problems related to entire transcendental functions) \cite{Speiser1929} (1929). He stresses on the importance of this problem in his paper \emph{Quelques propriétés des surfaces de Riemann correspondant aux fonctions méromorphes} (Some properties of Riemann surfaces corresponding to meromorphic functions) \cite{Ahlfors-JMPA} (1932), where he writes: ``This problem is, or ought to be, the central problem in the theory of functions. It is evident that its complete solution would give us, at the same time, all the theorems which have a purely qualitative character on meromorphic functions.''
   
 Nevanlinna,\index{Nevanlinna, Rolf (1895--1980)}  in the introduction of his book \emph{Analytic functions} (1953), writes (p. 1 of the English edition \cite{N-analytic}): ``[\ldots] Value distribution theory is thus integrated into the general theory of conformal mappings. From this point of view the central problem of the former theory is the \emph{type problem},\index{type problem} an interesting and complicated question, left open by the classical uniformization theory.''
  
  Teichm\"uller addressed the type problem\index{type problem} in several papers. We shall report here on his work on this subject in his two papers \emph{Eine Anwendung quasikonformer Abbildungen auf das Typenproblem} (An application of quasiconformal mappings to the type problem) \cite{T9}, published in 1937, and \emph{Untersuchungen \"uber konforme und quasikonforme Abbildungen} (Investigations on conformal and quasiconformal mappings) \cite{T13}, published in 1938.


 The Riemann surfaces that are considered by Teichm\"uller in these two papers are coverings of the sphere of the sphere branched over a finite number of values.
Nevanlinna,\index{Nevanlinna, Rolf (1895--1980)}  in his papers \cite{N-Ueber, N-Ueber-1932} and in  his book \cite{N-analytic} (in particular in  \S XI.2 and \S XII.1), studies a particular class of such surfaces. A particularly important class of examples is constituted by surfaces associated with automorphic functions, but there are many others. As a matter of fact, the idea of studying a Riemann surface which is a branched cover of the sphere and the problem of reconstructing it from the combinatorial information at the branch points and branch values originates with Riemann. The latter introduced Riemann surfaces in his doctoral dissertation \cite{Riemann-Grundlagen} and described, for the first time,  the domain of a meromorphic function to be a branched covering of the sphere. In this setting, he talked about surfaces ``spread over the sphere" (\emph{\"Uberlagerungsfl\"ache}). The problem of reconstructing the Riemann surface from combinatorial information at the branch points is a form of the so-called \emph{Riemann existence problem}\index{Riemann's existence problem} \cite{Riemann-Abelian}; see also Sto\"\i lov's work \cite{Stoilov} on the so-called \emph{Brouwer problem},\index{Brouwer problem} asking for a topological characterization of an analytic function. All this includes the kind of problems considered in Teichm\"uller's papers in a broader contex which can be traced back to the origin of Riemann surface theory.

      We shall also review some fundamental tools that appear in Teichm\"uller's  papers \cite{T9}  and \cite{T13}, in particular the  notion of line complex\index{line complex} that was extensively used by Nevanlinna\index{Nevanlinna, Rolf (1895--1980)}  and others, which is a topological/combinatorial tool that describes a type of branched covering.
       Independently of the work reviewed here, the ideas around the theory of line complexes\index{line complex} associated to Riemann surfaces should be interesting to low-dimensional geometers and topologists.
      
       Teichm\"uller, in his work on the type problem,\index{type problem} also relied on quasiconformal mappings. As a matter of fact, the paper \emph{Eine Anwendung quasikonformer Abbildungen auf das Typenproblem} \cite{T9}, whose subject is the type problem, is the first one written by Teichm\"uller in which he uses these mappings.  Teichm\"uller made extensive use of quasiconformal mappings in his approach to the type problem\index{type problem} as well as in other works he did on complex analysis, namely,  his papers on Riemann's moduli problem \cite{T20}   \cite{T29} \cite{T24}, on the Bieberbach coefficient problem \cite{T14} and on value distribution theory \cite{T13} \cite{T14}.

  \section{Line complexes}\label{s:line}
  
        We recall the definition of a  line complex\index{line complex} associated with a branched covering $\frak{M}$ of the Riemann sphere. This is a graph embedded in $\frak{M}$ which carries the combinatorics of this covering.   We shall follow the exposition in Nevanlinna's book \cite{N-analytic}, which is also the one Teichm\"uller uses in the two papers \cite{T9} and \cite{T13}.
  
   We assume that $\frak{M}$ is a simply connected Riemann surface which is a branched covering of the Riemann sphere $S^2=\mathbb{C}\cup\{\infty\}$, branched over finitely many points $a_1,\ldots,a_q\in S^2$ ($q\geq 2$) called the \emph{branching values}. Such a point is either of algebraic type\index{algebraic type (branching value)} (that is, the corresponding branching degree  is finite) or of logarithmic type\index{logarithmic type (branching value)}  (that is, the corresponding branching degree is infinite). Nevanlinna  considers the class of such surfaces in \cite[Chapter XI, \S 2 and \S 3]{N-analytic} and he denotes a surface belonging to this class by $F_q$, or by $F\left(a_1 , \cdots , a_q \right)$ if the branching values are prescribed. We shall keep using Teichm\"uller's notation $\frak{M}$. One example of such a surface is the universal covering of the sphere punctured at a finite set of points $a_1,\ldots,a_q$, $q\geq 2$ (see \cite[\S I.3]{N-analytic} for further details). Nevanlinna\index{Nevanlinna, Rolf (1895--1980)}  says that this surface is distinguished among the others as being ``as ramified as possible" \cite[p. 290]{N-analytic}.

    To define the line complex associated to $\frak{M}$, we start by drawing a simple closed curve $\gamma$ passing through the points $a_1,\ldots,a_q$ in the given cyclic order (chosen arbitrarily). This curve decomposes the sphere into two simply connected regions $\frak{I}$ and $\frak{A}$ that have a natural polygonal structures with vertices $a_1,\ldots,a_q$  and sides $(a_1a_2), (a_2a_3),\ldots (a_qa_1)$. The pre-image of $\gamma$ by the covering map $\frak{M}\to S^2$ is a graph $G$ which divides $\frak{M}$ into a finite or infinite number of pieces that are equipped with a natural polygonal structure induced from either $\frak{I}$ or $\frak{A}$. These pieces are called ``half sheets," and are labeled by $\frak{I},\frak{A}$.

     The line complex\index{line complex} is a graph embedded in the surface $\frak{M}$ which encodes the way the half sheets are glued together. We recall now more precisely the definition.
      
        If $\psi:   \frak{M}\to S^2$  denotes the branched covering map, then this map restricted to the set $\frak{M}\setminus \psi^{-1}(\{a_1,\ldots,a_q\})$ is an unbranched covering map onto its image $S^2\setminus \{a_1,\ldots,a_q\}$.
        A point lying above a branch value $a_i$ ($i=1,\ldots,q$) (or, equivalently, a vertex of a half sheet in $ \frak{M}$, allowing also a vertex to be at infinity), may be of three kinds:
            \begin{enumerate}
            \item Infinite order: such a vertex is not in $\frak{M}$ but lies at infinity.
            \item Finite order: there are $2m$ half sheets glued cyclically at such a point; in this case the point is said to be of order $m-1$.
            \item An unbranched point, called also ``nonliteral''. These are points of order zero.
            \end{enumerate}

    We construct a dual graph to the polygonal curve $\gamma$ in the sphere $S^2$ by choosing two points $P_1$ and $P_2$ in the interior of the polygons $\frak{I}$ and $\frak{A}$ respectively and joining, for each edge $(a_ia_{i+1})$ of $\gamma$, $P_1$ and $P_2$ by a simple arc that crosses $\gamma$  at a single point in the interior of that edge. This defines the dual graph of $\gamma$. Its lift by $\psi^{-1}$, as a planar graph embedded in $\frak{M}$, is the \emph{line complex}\index{line complex} $\Gamma$ of the covering.

 The line complex\index{line complex} $\Gamma$  might be finite or infinite 
 (in the sense that it may have a finite or infinite number of edges) depending on whether the Riemann surface $\frak{M}$ is a finite or infinite branched covering of the sphere. As a graph, it is homogeneous of degree $q$, and it is bipartite: choosing a coloring, say black and white, for the two vertices $P_1$ or $P_2$, the vertices of $\Gamma$ can also be colored black and white, depending on whether they are lifts of $P_1$ and $P_2$, and such that any edge of $\Gamma$ joins a white vertex to a black vertex. 
 
 The line complex $\Gamma$ decomposes the surface $\frak{M}$ into pieces called ``elementary regions" (the complementary components of the graph imbedded in that surface). Each elementary region is a polygon. 
  Such a polygon corresponds by 
 the covering projection $\psi$ to a unique branch value $a_i$ on the sphere. Corresponding to the above classification of the branch  points, the polygons are of three types:

    \begin{enumerate}
    \item A polygon with an infinite number of sides. It is unbounded in $\frak{M}$ and it corresponds to a logarithmic ramification point.
            \item A polygon with an even number $2m$ of sides, containing a ramification point of finite order $m-1$.
            \item A bigon (containing an unramified point over an $a_i$).
                        \end{enumerate}
        
   Line complexes\index{line complex} were introduced by  Nevanlinna\index{Nevanlinna, Rolf (1895--1980)}  in \cite{N-Ueber}  and Elfving in \cite{Elfving1934}\index{Elfving, Gustav (1908--1984)}. They were used as an important combinatorial  tool by Nevanlinna, \index{Nevanlinna, Rolf (1895--1980)} Ahlfors,\index{Ahlfors, Lars (1907--1996)} Speiser,\index{Speiser, Andreas (1885--1970)} Elflving,\index{Elfving, Gustav (1908--1984)} Ullrich,\index{Ullrich, Egon (1902--1957)} Teichm\"uller\index{Teichm\"uller, Oswald (1913--1943)} and other mathematicians in their investigation of the type problem.\index{type problem} In this setting, the type problem is reduced to that of recovering the type of a Riemann surface from the properties of its line complex\index{line complex}.
      
    As Nevanlinna and Ahlfors put it several times, the idea of encoding a branched covering of the sphere by a graph is contained in the paper \cite{Speiser1930} by Speiser.\index{Andreas Speiser (1885--1970)}\footnote{In the paper \cite{Ahlfors1982}, Ahlfors writes: ``Around 1930 Speiser had devised
a scheme to describe some fairly simple Riemann surfaces by means of a graph
and had written about it in his semiphilosophical style".} In fact, the latter introduced an object that was called later 
   \emph{Speiser tree}\index{Speiser tree}\index{Speiser tree}, as a topological tool to study the type problem.\index{type problem} The definition of a Speiser tree is close to that of a line complex, but we do not need to recall here the difference between the two notions.
      
 The branched covering $\frak{M}$ of the sphere is uniquely determined by the points $a_1,\ldots,a_q$, the simple closed curve $\gamma$ joining them and the line complex\index{line complex} $\Gamma$. In the papers \cite{Elfving1934} by Elfving\index{Elfving, Gustav (1908--1984)} and \cite{N-Ueber-1932} by Nevanlinna,\index{Nevanlinna, Rolf (1895--1980)}  Riemann  surfaces are constructed  from a given line complex satisfying certain conditions. There are sections on line complexes in the books by Nevanlinna \cite{N-analytic} and \cite{Nevanlinna&Paatero}, by Sario and Nakai \cite{SN}, by Goldberg and Ostrovskii \cite{goldberg}, and in other works. 

   \section{Teichm\"uller's paper \emph{Eine Anwendung quasikonformer Abbildungen auf das Typenproblem}} \label{s:Eine}

Teichm\"uller's motivation in his paper \emph{Eine Anwendung quasikonformer Abbildungen auf das Typenproblem} \cite{T9} is the study of the type problem\index{type problem} based on the properties of line complexes.\index{line complex} He  asks: ``How can one determine, from the given line complex, if the corresponding surface can be mapped one-to-one and conformally onto the whole plane, the punctured plane or the unit disc?"  He 
 notes that ``one is still very far away from sufficient and necessary criteria.''

 Let us start by making a list of the main ideas and results contained in that paper:

\begin{enumerate}
\item The question of determining the type of a Riemann surface from the properties of its line complex.\index{line complex}

\item \label{item2} The formula $$
D=|\mathcal K|+\sqrt{\mathcal K^2-1}\ \text{with} \ \mathcal K=\dfrac12\dfrac{u_x^2+u_y^2+v_x^2+v_y^2}{u_xv_y-v_xu_y}.
$$ for the \emph{dilatation quotient}\index{dilatation quotient} at a point of a map $x+iy\mapsto u+iv$ with continuous partial derivatives, and the definition of a quasiconformal mapping as a mapping with bounded dilatation quotient.\index{quasiconformal mapping}\index{dilatation quotient}

\item \label{i0}  The fact that the type of a simply connected Riemann surface is a quasiconformal  invariant.\index{type invariance under quasiconformal mapping}

\item \label{i1} The fact that the unit disc and the complex plane are not quasiconformally equivalent. 

\item \label{i2}  The fact that any two simply connected Riemann surfaces defined as branched coverings of the Riemann sphere with finitely many branch values and which have the same line complex\index{line complex}  
are quasiconformally equivalent.\index{quasiconformal mapping} 

\item \label{i3}  Techniques for piecing together quasiconformal mappings.

\item \label{i4}  The extension of an orientation-preserving diffeomorphism of the circle of class $C^1$ to a quasiconformal mapping of the disc which is \emph{conformal at the boundary}.  
 
\end{enumerate}

 These items are inter-related. For instance,  (\ref{i0}) is equivalent to (\ref{i1}),  (\ref{i4}) is used in the proof of (\ref{i3}), etc.

 We now discuss in more detail the content of the paper \cite{T9}.

Teichm\"uller proves the following:
\begin{theorem}
 Any two Riemann surfaces with the same line complex\index{line complex} have the same type.
\end{theorem}

In order to prove this result, Teichm\"uller shows the following key results which we state as propositions:
\begin{proposition}\label{th:same}
Any two branched coverings of the sphere with equal\index{line complex} line complexes can be mapped quasiconformally onto each other. 
\end{proposition}

\begin{proposition}\label{th:invariant}
 The type of a simply connected open Riemann surface is invariant under quasiconformal mappings.
\end{proposition}

 The latter proposition is a direct consequence of the following result:
\begin{proposition} \label{th:cannot}
 The complex plane and the unit disc cannot be conformally mapped onto each other.
\end{proposition}

After stating his results, Teichm\"uller, in §2 of his paper, defines the notion of ``dilatation quotient"\index{dilatation quotient} of a mapping at a point where the mapping is of class $C^1$.

We recall that the differential of a $C^1$ mapping takes an infinitesimal ellipse centered at a point to an infinitesimal ellipse centered at the image.  The
dilatation quotient of a $C^1$  mapping $f$ is a function that assigns to each point of the domain of $f$ the ratio of the major axis to the minor axis of an ellipse that is the image of a circle  centered at the origin by the differential of $f$. This dilatation quotient is also known as the \emph{distortion of the Tissot indicatrix};\index{Tissot indicatrix}  see Gr\"otzsch's\index{Gr\"otzsch, Herbert (1902--1993)} paper  \cite{Gr1930} where the author uses this expression, and the review \cite{Pa-Tissot} on Tissot's work in this volume. There is also a relation with the \emph{characteristic functions}\index{characteristic (Lavrentieff)}  introduced by Lavrentieff\index{Lavrentieff, Mikha\"\i l (1900-1980)} in \cite{Lavrentiev1935a} and \cite{Lavrentiev1935}. Teichm\"uller derives formula (\ref{item2}) above,  which allows one to compute the dilatation quotient in terms of the partial derivatives of the mapping.  

In §3, Teichm\"uller first defines a quasiconformal mapping\index{quasiconformal mapping} as a continuous one-to-one and sense-preserving mapping between two domains of the plane that is differentiable except at isolated points and whose dilatation quotient\index{dilatation quotient}  is bounded. This definition is close  to the one given by Gr\"otzsch in \cite{Gr1928}. The differentiablity conditions in the definition of a quasiconformal mapping were relaxed later; cf. for instance Lehto's historical survey on quasiconformal mappings \cite{lehto}.  

 Teichmüller notes that the definition of a quasiconformal mapping between domains of the plane can be transferred to the setting of mappings between Riemann surfaces. He then gives conditions for piecing together continuously differentiable and quasiconformal  mappings defined on  two domains $\mathfrak{S}_1$ and $\mathfrak{S}_2$ separated by a curve $\mathfrak{C}$. He concludes this section  by solving what he describes as ``an important special mapping problem,'' namely, the problem of extending an arbitrary orientation-preserving $C^1$ diffeomorphism of the unit circle $\mathbb{S}^1$ to a quasiconformal mapping\index{quasiconformal mapping} of the unit disc $\mathbb{D}$ that is ``conformal at the boundary.''\index{quasiconformal mapping!conformal at the boundary} Here, a differentiable mapping of the unit disc is said to be conformal at the boundary if its dilatation quotient on sequences of points approaching the boundary approaches $1$.

\begin{remarks} 1. The extension that Teichm\"uller constructs is (up to a clockwise rotation by $\frac{\pi}{4}$)  the one that appears  in Ibragimov's paper \cite{Ib2010}, in which the latter says that  such an extension was known to Ahlfors and Beurling (without any further reference). Presumably, Ibragimov was not aware of Teichm\"uller's result.

2.  A result of Ahlfors and Beurling says that a circle homeomorphism can be extended quasiconformally to the unit disc if and only if it is quasisymmetric (see \cite{BA1956}). Today, two well-known quasiconformal extensions\index{quasiconformal extension} of quasi-symmetric homeomorphisms of the circle are the \emph{Beurling--Ahlfors}\index{Beurling--Ahlfors extension} and the \emph{Douady--Earle}\index{Douady--Earle extension} extensions, introduced in the papers \cite{BA1956} and \cite{douady&earle} respectively. The two extensions  are different from the one described by Teichm\"uller. 

3. Teichm\"uller's notion of quasiconformal mappings of the disc that are ``conformal at the boundary" is used today in the definition of ``asymptotically conformal"\index{asymptotically conformal quasiconformal mapping}\index{quasiconformal mapping!asymptotically conformal} quasiconformal mappings of the disc, a class of mappings that are at the basis of the notion of \emph{asymptotic Teichm\"uller space of the disc} introduced by Gardiner and Sullivan in their paper \cite{GS}.
\end{remarks}
 
 In §4, Teichm\"uller proves  that  two surfaces having the same line complexes\index{line complex} can be mapped quasiconformally  onto each other (Proposition \ref{th:same} above). With the notation introduced for line complexes in \S \ref{s:line}, after mapping the two discs $\frak{A}$ and $\frak{I}$ quasiconformally onto the inside and the outside of the unit disc respectively, he uses the extension he introduced in the previous section in such a way that the mapping obtained is \emph{conformal on the simple closed curve $\gamma$} (the common boundary of the discs $\frak{A}$ and $\frak{I}$), and defining therefore a quasiconformal mapping between the $z$-plane and the $w$-plane. He performs the same operation with another Riemann surface $\frak{M}^\prime$ that has the same line complex as $\frak{M}$. In this way, he gets a quasiconformal mapping from a $z^{\prime}$-plane onto the $w$-plane and then a quasiconformal mapping from the $z$-plane onto the $z^{\prime}$-plane. Because $\frak{M}$ and $\frak{M}^\prime$ share the same line complex,\index{line complex} the gluing specification between corresponding half sheets is the same.  Teichm\"uller eventually gets a quasiconformal mapping between $\frak{M}$ and  $\frak{M}^{\prime}$.

 In §5, he proves that the unit disc and the complex plane cannot be quasiconformally mapped onto each other (Proposition \ref{th:cannot} above). The proof is by contradiction and it is based on the so-called  \emph{length-area} method,\index{length-area method} a method to which Teichm\"uller referred in his later papers as the \emph{Gr\"otzsch--Ahlfors}\index{Gr\"otzsch--Ahlfors method}  method (see for instance \cite{T13, T20, T23,T24}). More precisely, Teichm\"uller reasons by contradiction. He assumes that there exists a quasiconformal mapping from the complex plane to the unit disc. Such a mapping sends, for any $0<r_1 < r_2$, the annulus $\left\lbrace z \mid r_1 < \left| z \right| < r_2 \right\rbrace$ onto  a doubly-connected subdomain of the unit disc.  The image being bounded, it has a bounded conformal modulus and therefore the modulus of  $\left\lbrace z \mid r_1 < \left| z \right| < r_2 \right\rbrace$  is also bounded. This is obviously false since $r_2$ can be chosen arbitrarily large. 
 
In an amendment made during the page proofs of his paper and written at the end of §5, Teichm\"uller notes that the final part of the paper \cite{Gr1928} by Gr\"otzsch\index{Gr\"otzsch, Herbert (1902--1993)} published in 1928 is closely related to his own proof of Proposition \ref{th:invariant} (information he obtained from H. Wittich).\index{Wittich, Hans}\footnote{A translation of Gr\"otzsch's paper is provided in the present volume.} In that paper, Gr\"otzsch uses this method in order to extend the small and the big Picard theorems to quasiconformal (and not only conformal) mappings. Lavrentieff\index{Lavrentieff, Mikha\"\i l (1900-1980)} uses the same idea in his 1935 paper \cite{Lavrentiev1935} which we mention below.  Conceivably,  Teichm\"uller was not aware of Lavrentieff's paper.\footnote{Teichm\"uller was aware of Lavrentieff's paper when he wrote the later paper \cite{T13}; see the footnote in §6.4 of \cite{T13}.}

 In the last section, Teichm\"uller provides an example to briefly show that for a simply connected Riemann surface $\frak{M}$ defined as a branched covering of the sphere with branch values $a_1 , \cdots a_q$, the position of these points and the simple closed curve $\gamma$ joining them that is used in the construction of the line complex, although it has no influence on the question of deciding the type of the surface, has an influence on the value distribution of the uniformizing mapping, namely, the biholomorphic mapping that sends $\frak{M}$  onto the disc or to the complex plane.\index{Nevanlinna theory}\index{value distribution theory} He deals with the explicit example of a branched covering with four branch values and with a line complex determined by the intersection of the real line with infinitely many circles. He gives explicitly the \emph{Nevanlinna characteristic}\index{Nevanlinna!characteristic} of the associated meromorphic mapping and shows that the coefficients of this characteristic depend on the  choice of $\gamma$.  He declares that the proof is not difficult but will be carried out in a different context. The same example is further considered in his 1944 paper \cite{T33}.

  Teichm\"uller was not the first to be interested in the application of quasiconformal mappings to the type problem,\index{type problem} and we refer here again to the 1935 work of Lavrentieff\index{Lavrentieff, Mikha\"\i l (1900-1980)} in \cite{Lavrentiev1935} (see the English translation in the present volume; cf. also \cite{Lavrentiev1935a} for a summary of the results and the paper \cite{AP-Lavrentiev}  for a commentary also in the present volume). Lavrentieff\index{Lavrentieff, Mikha\"\i l (1900-1980)} introduced the notion of \emph{fonction presque analytique} (almost analytic function),\index{almost analytic function} a class of mappings that are more general than the ones introduced earlier by Gr\"otzsch,\index{Gr\"otzsch, Herbert (1902--1993)} and whose dilatation quotient\index{dilatation quotient} is not necessary bounded. He gives, for a particular class of Riemann surfaces defined as the graph of a function of two variables, a sufficient condition on the dilatation quotient that makes the surface of hyperbolic type (see \cite[Th\'eor\`{e}me 8]{Lavrentiev1935}). The same year Teichm\"uller  published his paper \cite{T9} (1937), Kakutani published a paper \cite{Kakutani} in which he also studies the type problem\index{type problem} using (another) notion of quasiconformal mappings. See also Kobayashi's paper \cite{Kobayashi}.

In his paper \cite{T13} which we consider in \S \ref{s:Unter}, Teichm\"uller gives  a criterion to find the type of Riemann surfaces in a class already considered in his paper \cite{T9}, and as a consequence he constructs an example of a hyperbolic Riemann surface with a so-called ``mean ramification" equal to $2$. We review this in the next section.

\section{Nevanlinna's conjecture}\label{s:Nevanlinna}
 
           The notation is the one introduced in \S \ref{s:line}. Nevanlinna's conjecture is based on the fact that the type of the simply connected surface $\frak{M}$ may be deduced from information on the measure of the amount of ramification of its associated line complex.\index{line complex!degree of ramification} Roughly speaking, if the amount of ramification is small (in some sense to be defined), the surface is of parabolic type, and if it is large, it is of hyperbolic type. The motivation is that a large degree of ramification will put a large angle structure at the vertices of the graph $G$. This can be made precise by introducing a notion of combinatorial curvature, and it is the object of these investigations. Nevanlinna\index{Nevanlinna conjecture}\index{conjecture!Nevanlinna} writes that ``it is natural to imagine the existence of a critical degree of ramification that separates the more weakly branched parabolic surfaces from the more strongly branched hyperbolic surfaces"  \cite[p. 308]{N-analytic}. To make precise this observation, he introduces a notion of curvature of a line complex, which he calls the \emph{mean excess}.\index{line complex!mean excess} His conjecture says that the surface $\frak{M}$ is parabolic or hyperbolic depending on whether the mean excess of its line complex is zero or negative respectively.  Teichm\"ulller in his paper \cite{T13}, disproves the conjecture by exhibiting a hyperbolic simply connected Riemann surface branched over the sphere 
 whose mean excess of the corresponding line complex\index{line complex} is equal to zero. 
                    
                  Nevanlinna's study of this topic is restricted to the class of surfaces  $F_q$ that we defined earlier (simply connected surfaces that are branched covers of the sphere and whose branch points project to a finite number of points $a_1,\ldots,a_q$). We quote him from his book  \cite[p. 312]{N-analytic}:
   \begin{quote}\small
   [\ldots] Relative to transcendental surfaces $F_q$ the following question arises:
   \emph{Is the surface parabolic or hyperbolic according to the angle geometry of the surface $F_q$ is ``euclidean" or ``Lobachevskyan," i.e. according as the mean excess $E=2-V$ is zero or negative?} This question was answered (1938) by Teichm\"uller \cite{T13} and, indeed, in the negative.
   \end{quote}
   
    In the rest of this section, we explain Nevanlinna's conjecture\index{Nevanlinna conjecture}\index{conjecture!Nevanlinna}.

                     We use the notation for line complexes introduced in \S \ref{s:line}. The vertices of the complex $\Gamma$ correspond to the half sheets in $\frak M$ and the faces of $\Gamma$, which are all homeomorphic to discs, are called \emph{elementary polygons}. The number of sides of an elementary polygon may be infinite, but if it is finite, then it is even. Indeed, an elementary polygon covers a cell in the Riemann sphere whose boundary consists of two edges. Therefore, each face has $2m$ edges, with $m=1,2,\ldots,\infty$. Each elementary polygon with $2m$ sides is associated with a branch point of order $m-1$.

In order to give an idea of the ramification measure that he introduces, Nevanlinna first considers the case of a closed surface which is a finite-sheeted covering of the sphere. In this case, the ramification measure is the sum of the orders of the branch points divided by the number of sheets of the covering. The calculation in this special case is made in terms of the associated line complex.\index{line complex} 
                                                
                                               The goal is to assign a ramification measure to each vertex of the graph. (Recall that the ramification points of the covering map $\psi$ from  $\mathfrak{M}$ to the Riemann sphere are at the ``centers" of the polygons that are the complementary components of the line complex.)
                                                
                                                To each branch point of order $m-1$ corresponds an elementary polygon with $2m$ sides and $2m$ vertices. Twice the order of the branch point, that is, $2m-2$, is distributed onto the $2m$ vertices of the polygon. In this way, each vertex $P$ of $\Gamma$ receives from the given elementary polygon a contribution of the line complex $\frac{2m-2}{2m}=1-\frac{1}{m}$ of the order of ramification. Therefore, the total ramification at a vertex $P$ is equal to
                                                \[V_P=\sum_W(1-\frac{1}{m}),\]
                                                the sum, for each  $P$, being taken over the  elementary polygon $W$ having $P$ as a vertex. 
                                                
                                                Since the number of vertices of $\Gamma$ is equal to $2n$ (twice the number of sheets of the covering), the mean ramification of the surface, which is naturally defined as
                                                \[\frac{1}{2n}\sum_PV_P,\]
                                                is equal to 
                                                  $$
                                                  2-\frac{2}{n}.
                                                  $$

                        We now return to our simply connected infinitely-sheeted surface $\frak{M}$, branched over $q$ points, with associated line complex $\Gamma$. 

As in the previous special case, to each vertex $P$ of $\Gamma$ is associated a ramification index $V_P$ coming from its adjacent elementary polygons. With the above notation, the index coming from one polygon is equal to $1-\frac{1}{m}$ and since the number of edges of the elementary polygons may be infinite, we allow $m$ to be infinite, in which case the ramification index coming from this region is 1.

Since $P$ is a vertex of at least 2 and at most $q$ adjacent elementary regions whose order $m-1$ is positive, we always have  $1\leq V_P\leq q$.

Since the surface $\frak M$ has infinitely many sheets, to define the average ramification of the associated line complex $\Gamma$, Nevanlinna uses what he calls a \emph{wreathlike exhaustion} (see \cite[\S XII.1.3]{N-analytic}) of this graph by an infinite sequence of subgraphs $\Gamma^\nu$ ($\nu=1,2,\ldots$). This is done by starting with a vertex $P_0$ chosen arbitrarily, which we call the base vertex, then defining the complex $\Gamma^1$ as the subgraph of $\Gamma$ obtained by adding to $P_0$ the vertices $P_1$ at distance one from this base vertex (first generation) together with the edges joining these new vertices to $P_0$, then the complex $\Gamma^2$ as the one obtained by adding the vertices $P_2$ at distance one from $\Gamma^1$ (second generation) together with the edges joining these new vertices to $\Gamma^1$, etc. Here, the distance between two vertices is defined as the minimal number of edges joining them.

For every $\nu=1,2,\ldots$, we let $n_\nu$ be the number of vertices in the approximating graph $\Gamma^\nu$. The \emph{mean ramification}\index{mean ramification} of $\Gamma^\nu$  is defined as
\[V_\nu=\frac{1}{n_\nu}\sum_{F^\nu}V_P.\]

 The mean ramification of $\Gamma$ (and of the surface $\frak{M}$) is the limit
\[V=\lim_{\nu\to\infty}V_\nu,\] provided this limit exists. If this limit does not exist, Nevanlinna works with the lower limit $\varliminf V_\nu$ and the upper limit $\varlimsup V_\nu$.

There is a geometric way of calculating  this measure of ramification, in relation with a combinatorial notion of curvature. We return to the graph $G=\psi^{-1}(\gamma)$ in $\frak{M}$, that is, the pre-image by the projection map of the simple closed curve $\gamma$ in the Riemann sphere that joins the branchvalues $a_1,\ldots,a_q$. We recall that $G$  defines a combinatorial decomposition of $\frak{M}$ whose cells have $q$ sides  and whose vertices lie above the branching values $a_1,\ldots a_q$. There is a natural combinatorial angle structure associated to this polygonal decomposition of $\mathfrak{M}$, where to each vertex of $G$ which corresponds to a branch point of order $m-1$ one associates an angle measure equal to $\pi/m$. The \emph{excess}\index{excess} of a polygon $P$ with $q$ sides is then defined as 
\[E_P= \sum \frac{1}{m}-q+2\]
where the sum is taken over the angles of the polygon. 
This is justified as follows:

If the $q$-sided polygon were a Euclidean polygon, then its angle sum would be equal to $(q-2)\pi$. Suppose now that we are given a $q$-sided polygon which is equipped with a metric or more generally with a structure that assigns to it a reasonable angle structure at the vertices. Then, the ``excess" of this polygon, which is a measure of its deviation from being Euclidean, is taken to be its angle sum divided by $\pi$ minus the angle sum of a Euclidean polygon with the same number of sides divided by $\pi$. This is consistent with the formula we gave for $E_P$.

The above formula for the ramification $V_P$ of a polygon $P$ gives
\[V_P=\sum(1-\frac{1}{m})=q-\sum\frac{1}{m}.\]
Thus, we have
\[V_P+E_P=2,\]
that is, the sum of the ramification and the excess of a fundamental polygon $P$ is equal to 2.

To get a more precise idea of the geometric meaning of ramification, we consider the case of a regularly ramified surface, where $V=V_P$ for every fundamental polygon $P$. There are three cases:
\begin{enumerate}
\item Elliptic case: $V<2$. One can take the fundamental polygons on $\frak{M}$ to be spherical polygons with geodesic sides and angle excess $\pi E$ ($E=E_P$ being the excess defined above).

\item Parabolic case:  $V=2$, $E=0$. The fundamental polygons can be taken as Euclidean polygons, with zero excess.

\item Hyperbolic case: $V>2$, $E<0$. The fundamental polygons can be taken as hyperbolic polygone, with angle deficit equal to $\pi\vert E\vert$.

\end{enumerate}

This motivates the following conjecture which Nevanlinna made in the general case of a non-necessarily regular surface \cite[p.312]{N-analytic}:
The surface is parabolic\index{parabolic surface} (respectively hyperbolic)\index{hyperbolic surface} if the mean excess\index{line complex!mean excess} $E=2-V$ is zero (respectively negative).

                       Teichm\"uller, in \S 7 of the paper \cite{T13} disproves the conjecture by exhibiting a hyperbolic simply connected surface $S$ branched over the sphere and whose  mean excess is zero. 
                       We review Teichm\"uller's work in the next section.

                       Let us note that in a more recent paper \cite{BMS}, the authors disprove the second part of the conjecture by giving an example of a simply connected Riemann surface which is parabolic and which has negative mean excess.

    \section{The type problem in Teichm\"uller's paper \emph{Untersuchungen \"uber konforme und quasikonforme Abbildungen}}  \label{s:Unter}
   
In this section, we review the criterion for hyperbolicity of surfaces in a certain class given by  Teichm\"uller  in the last section of his paper \emph{Untersuchungen \"uber konforme und quasikonforme Abbildungen}. The section is titled \emph{A type criterion}.

                       Teichm\"uller considers a class of surfaces  $\frak{M}$  satisfying the following properties (we use the notation of the previous section):
                       
                       \begin{enumerate}
                       \item no branch point is algebraic;
                       \item the line complex\index{line complex} $\Gamma$ associated to $\frak{M}$ does not have any infinite chain of edges with no branch points;
                       \item each vertex of $\Gamma$ has either two or $q$ edges starting from it.

                                      \end{enumerate}
                                             
                      Teichm\"uller proves that such a surface is hyperbolic if and only if it is sufficiently ramified, and for this, introduces a notion of measure of ramification\index{measure of ramification} which is different from the one introduced by Nevanlinna (\S \ref{s:Nevanlinna} above). We recall the definition. 
                      
                      We take a vertex $B_0$ adjacent to $q$ simple edges (by Conditions (2) and (3), such a vertex exists). From $B_0$, we follow, from each edge adjacent to it, a chain of edges until we arrive at a branch vertex. We define in this way $q$ vertices that are branch points, and we call this set of vertices $B_1$. By Condition (3), each vertex in $B_1$ has $q-1$ outgoing branches. Continuing in the same manner, we obtain $q(q-1)$ branch vertices $B_2$, then $q(q-1)^2$ branch vertices $B_3$, etc. Thus, for each $k\geq 1$, there are $q(q-1)^{k-1}$ \emph{chains} leading from the set of vertices $B_{k-1}$ to the set of vertices $B_k$. Condition (1) says that continuing in this way we  traverse the entire line complex without repetitions.

                       Each chain of vertices has a well-defined length, and Teichm\"uller denotes by $\varphi(k)$ the maximum of the lengths of all chains connecting an element in $B_k$ to an element in $B_{k-1}$. For each $k\geq 1$, he sets
                      \[\psi(k)=\Max_{\kappa=1,\ldots,k}\varphi(k).\]
                      
                      The smaller $\varphi(k)$ and $\psi(\kappa)$ are, the more ramified is the surface $\frak{M}$.
                      
                      Teichm\"uller proves the following:
                      \begin{theorem}
                      If the series $\sum_{k=1}^\infty\frac{\psi(k)}{k^2}$ converges, then $\frak{M}$ is of hyperbolic type.
                      \end{theorem}
                                            
                      The  proof uses an explicit conformal representation onto the unit disc and hyperbolic geometry. It is based on the construction of a quasiconformal mapping (with a non-necessarily bounded dilatation) from $\mathfrak{M}$ onto the modular surface (the universal covering of the Riemann sphere punctured at the points $a_1,a_2,\ldots, a_q$), which is known to be hyperbolic. In order to conclude with the invariance of the type, Teichm\"uller uses a result he previously obtained (see \cite[§6.5]{T13}) and which is a generalisation of Proposition \ref{th:cannot}. The quasiconformal mapping from $\mathfrak{M}$ onto the modular surface is constructed by first modifying the line complex associated to $\mathfrak{M}$ and making it a \emph{Speiser tree},\index{Speiser tree} by replacing all chains without branching by bundles of edges arranged in a row.

                      Teichm\"uller gives then an example of a surface $\frak{M}$ whose mean ramification $V$ in the sense of Nevanlinna (\S \ref{s:line} above) is equal to 2, while $\sum_{k=1}^\infty\frac{\psi(k)}{k^2}$ converges. Thus, the surface is hyperbolic, which disproves Nevanlinna's conjecture.
                       
 In an addendum to his paper written on the page proofs, Teichm\"uller notes that Kakutani, in his paper \cite{Kakutani} (1937), obtained a sufficient criterion for hyperbolicity that applies to a class of  surfaces he considers  in \S \ 7.3 of his paper \cite{T13}. He says that Kakutani's criterion does not contain his own criterion, but gives a better order of magnitude. He then writes about Kakutani: ``If he had worked with unbounded dilatation coefficients and applied my \S \ 6.5, he would have obtained a much better criterion for \S \ 7.3 than the one developed here. His quasiconformal mapping is hence more suitable than mine; it is (for now) restricted to the special surfaces in  \S  7.3."

 \begin{remarks}
         
         1.-- The results on the type problem obtained in Teichm\"uller's papers  considered here concern surfaces which are branched coverings of the Riemann sphere with logarithmic branch points. There is a result of Ahlfors on the type problem for branched coverings with only algebraic branch points, see \cite{Ahlfors1931} (1931) and results by Speiser \cite{Speiser1930, Speiser1932}\index{Speiser, Andreas (1885--1970)} and Ullrich \cite{U-S}\index{Ullrich, Egon (1902--1957)}.  Ahlfors'\index{Ahlfors, Lars (1907--1996)} first works on the type problem (see \cite{Ahlfors1931} and \cite{Ahlfors-CRAS}), unlike the works of Laverentieff on the same problem, involve the so-called length-area method which is a common tool used by Gr\"otzsch and himself, and later by Teichm\"uller, in relation to quasiconformal mappings.

2.-- We already noted that Lavrentieff, before Teichm\"uller, used quasiconformal maps (in the version he called ``almost analytic") in the study of the type problem (cf. \cite{Lavrentiev1935} and \cite{Lavrentiev1935a}). His approach was Riemannian, unlike the combinatorial approach of Teichm\"uller. See the translation of Lavrentieff's paper \cite{Lavrentiev1935} in the present volume and the corresponding commentary \cite{AP-Lavrentiev}).

4.-- There is a large literature on the type problem. Kakutani made the relation with Brownian motion. In \cite{Kaku} (1945), he proved the following:  Let $S$ be a simply-connected open Riemann surface which is an infinite cover of the sphere. Let $D$ be a simply connected subsurface of $S$ which is bounded by a Jordan domain $\Gamma$. For a point $\zeta$ in $S\setminus \overline{D}$, let $u(\zeta)$ be the probability that the Brownian motion on $S$ starting on $\zeta$ enters into $\Gamma$ without getting out of $S-\overline{D}$ before.  Then, one of the following two cases holds:
\begin{enumerate}
\item $u(\zeta)$ is $<1$ everuwhere in $S-\overline{D}$, and in this case $S$ is hyperbolic

\item $u(\zeta)$ is identically equal to $1$ on $S-\overline{D}$, and in this case $S$ is parabolic.
\end{enumerate}

 Doyle made a relation between the type problem, Brownian motion and the propagation of an electric current on the surface, see \cite{DS}. In fact, this is a return to Riemann's point of view. Indeed, the latter relied, at several places in his work on Riemann surfaces, to arguments from electricity (see Riemann's doctoral dissertation \cite{Riemann-Grundlagen} and his paper on Abelian functions \cite{Riemann-Abelian}. The interested reader may also refer to the survey \cite{Riemann1} on physics in Riemann's mathematical work).

  \end{remarks}

\bigskip
\bigskip

 \noindent {\emph{Acknowledgements.}}  The first and third authors acknowledge support from the U.S. National Science Foundation grants DMS 1107452, 1107263, 1107367 ``RNMS: GEometric structures And Representation varieties'' (the GEAR Network).

\bigskip

\bigskip

\textsc{Vincent Alberge, Fordham University, Department of Mathematics, 441 East Fordham Road, Bronx, NY 10458, USA}

\textit{E-mail address}: {\tt valberge@fordham.edu}

\textsc{Melkana Brakalova, Fordham University, Department of Mathematics, 441 East Fordham Road, Bronx, NY 10458, USA}

\textit{E-mail address}: {\tt brakalova@fordham.edu}

\textsc{Athanase Papadopoulos, Institut de Recherche Mathématique Avancée, Universit\'{e} de Strasbourg et CNRS, 7 rue Ren\'{e} Descartes, 67084 Strasbourg Cedex, France}

\textit{E-mail address}: {\tt papadop@math.unistra.fr} 
\printindex

\end{document}